\documentclass[12pt, a4paper, reqno]{amsart}

\usepackage[english]{babel}

\usepackage[T1]{fontenc}
\usepackage{lmodern}

\usepackage{amsmath}
\usepackage{amssymb}
\usepackage{amsthm}

\usepackage{url}

\usepackage{xcolor}
\usepackage{comment}
\usepackage{tikz-cd}

\usepackage{hyperref} 

\usepackage{enumitem}
\usepackage[normalem]{ulem}

\newcommand{\lexp}[2]{\null^{#2} \mkern-2mu #1}
\newcommand{\lexpp}[2]{\null^{#2} \mkern-2mu (#1)}

\renewcommand{\epsilon}{\varepsilon}
\renewcommand{\theta}[0]{\vartheta}
\renewcommand{\phi}[0]{\varphi}

\newcommand{\lift}[1]{#1^{\uparrow}}

\newcommand{\ZZ}{\mathbb{Z}}

\newcommand{\ppar}{\ \par}

\newcommand{\Span}[1]{\left\langle\, #1 \,\right\rangle}

\newcommand{\nothing}[1]{#1}

\DeclareMathOperator{\GL}{GL}

\DeclareMathOperator{\Aut}{Aut}

\DeclareMathOperator{\Inn}{Inn}

\DeclareMathOperator{\coh}{H}

\newtheorem{dummy}{Dummy}
\numberwithin{dummy}{section}
\numberwithin{figure}{section}

\newtheorem{theorem}[dummy]{Theorem}

\newtheorem{lemma}[dummy]{Lemma}

\newtheorem{proposition}[dummy]{Proposition}
\newtheorem{prop}[dummy]{Proposition}

\theoremstyle{definition}

\theoremstyle{remark}

\newtheorem{remark}[dummy]{Remark}

\newtheorem{question}[dummy]{Question}

\makeatletter
\def\imod#1{\allowbreak\mkern10mu({\operator@font mod}\,\,#1)}
\makeatother

\numberwithin{equation}{section}

\usepackage{amsthm, mathtools, amssymb, amsfonts, stmaryrd, latexsym,
mathrsfs, todonotes, color}

\begin{document}

\date{29 April 2022, 14:09   CEST --- Version 4.11%
}

\title[Skew braces and Rota--Baxter operators]
      {Skew braces from Rota--Baxter operators:\\
        A cohomological characterisation and\\
        some examples} 
      
\author{A.~Caranti}

\address[A.~Caranti]%
 {Dipartimento di Matematica\\
  Universit\`a degli Studi di Trento\\
  via Sommarive 14\\
  I-38123 Trento\\
  Italy\\\endgraf
  ORCiD: 0000-0002-5746-9294} 

\email{andrea.caranti@unitn.it} 

\urladdr{https://caranti.maths.unitn.it/}

\author{L.~Stefanello}

\address[L.~Stefanello]
        {Dipartimento di Matematica\\
          Universit\`a di Pisa\\
          Largo Bruno Pontecorvo, 5\\
          56127 Pisa\\
          Italy\\\endgraf
          ORCiD: 0000-0003-3536-2397}
        
\email{lorenzo.stefanello@phd.unipi.it}

\urladdr{https://people.dm.unipi.it/stefanello}

\subjclass[2010] {17B38 20J06 20E22 20D15 08A35}

\keywords{endomorphisms, skew braces, Rota--Baxter operators, group
  cohomology, central extensions, liftings}

\begin{abstract}
  Rota--Baxter operators  for groups  were recently  introduced by
  L.~Guo, H.~Lang,  and Y.~Sheng.

  V.~G.~Bardakov  and V.~Gubarev  showed that  with each  Rota--Baxter
  operator one can associate a skew  brace. Skew braces on a group $G$
  can be characterised in terms of certain gamma functions from $G$ to
  its automorphism group  $\Aut(G)$, that are defined  by a functional
  equation.  For the skew braces obtained from a Rota--Baxter operator
  the corresponding gamma functions   take values in the inner
  automorphism group $\Inn(G)$ of $G$.

  In this paper,  we give a characterisation of the
  gamma functions on a group $G$, with values in $\Inn(G)$, that come
  from a Rota--Baxter operator, in terms of the vanishing of a certain
  element in a suitable second cohomology group.

  Exploiting this characterisation, we are  able to exhibit examples
  of skew braces whose corresponding gamma functions take
  values in the  inner automorphism group, but cannot be obtained
  from a Rota--Baxter operator.

  For  gamma  functions  that  can be  obtained  from  a  Rota--Baxter
  operators,  we  show  how  to  get the  latter  from  the    former,
  exploiting the  knowledge that a suitable  central group extension
  splits.
\end{abstract}

\thanks{%
  Both authors are members of INdAM---GNSAGA. 
  The first  author gratefully acknowledges support  from the Department
  of  Mathematics  of  the  University of  Trento.   The  second  author
  gratefully acknowledges support from  the Department of Mathematics of
  the University of Pisa.}

\maketitle

\thispagestyle{empty}

\section{Introduction}

Rota--Baxter operators for various kinds of algebras have been studied
by several authors since G.~Baxter introduced them for commutative
algebras~\cite{Baxter} in 1960.

Recently, L.~Guo, H.~Lang, and Y.~Sheng introduced Rota--Baxter
operators for groups~\cite{GLS}. These were studied further by
V.~G.~Bardakov and V.~Gubarev~\cite{BarGub1, BarGub2}. In particular,
in~\cite{BarGub2}, it is  showed how to associate a skew brace with
a Rota--Baxter operator.

Recall that a  \emph{skew (left) brace}, defined  in~\cite{skew}, is a
triple $(G, \cdot, \circ)$, where $(G, \cdot)$ and $(G, \circ)$ are
groups and the two operations are related by the identity
\begin{equation*}
  g\circ(h\cdot k)
  =
  (g\circ h)\cdot g^{-1}\cdot (g\circ k)
\end{equation*}
for all $g,h,k\in G$; here the inverse refers to the operation
``$\cdot$''. 

The interplay between set-theoretic solutions of the Yang--Baxter
equation of Mathematical Physics, (skew) braces, regular subgroups,
and Hopf--Galois structures has spawned a considerable body of
literature in recent years; see, for example,~\cite{braces, skew, CJO,
  GP, Byo96, SV2018}.

Given a group $(G, \cdot)$, the group operations ``$\circ$'' such that
$(G, \cdot, \circ)$ is a skew brace can be characterised in the form
\begin{equation}
  \label{eq:circ}
  g \circ h = g\cdot \lexp{h}{\gamma(g)},
\end{equation}
where $\gamma \colon  G \to \Aut(G)$ is a function,  which we write as
an exponent, satisfying the identity
\begin{equation*}
  \gamma(g\cdot \lexp{h}{\gamma(g)})
  =
  \gamma(g) \gamma(h)
\end{equation*}
for all $g,h\in G$.
We call such a function a \emph{gamma function}. These functions
are usually referred to as $\lambda$ in the literature of skew
braces. (See for instance~\cite{p2q, CS1} for this context.)

A \emph{Rota--Baxter operator} on the group $(G, \cdot)$ is a function
$B \colon G \to 
G$ that  satisfies the identity
\begin{equation*}
  B(g\cdot B(g)\cdot h\cdot B(g)^{-1}) = B(g) \cdot B(h)
\end{equation*}
for all $g, h \in G$. 

Consider the morphism from the group $G$ onto the
group of its inner automorphisms
\begin{align*}
  \iota \colon\ G &\to \Inn(G)\\
  g &\mapsto (x \mapsto g\cdot x\cdot g^{-1}),
\end{align*}
whose kernel is the centre $Z(G)$ of $G$.
If $B$ is a Rota--Baxter operator on the group
$G$, then the function
\begin{equation*}
  \gamma(g) = \iota(B(g))
\end{equation*}
is immediately  seen to  be a  gamma function on  $G$, with  values in
the inner automorphisms group $\Inn(G)$. This function yields the skew
brace $G(B)$ introduced 
in~\cite{BarGub2}, where
\begin{equation*}
  g \circ h=g\cdot \lexp{h}{\iota(B(g))}
  =
  g\cdot B(g)\cdot h\cdot B(g)^{-1}
\end{equation*}
for all $g, h \in G$.

We address  here the question whether  the
converse holds.
\begin{question}
  \label{q:titq}
  Let $\gamma$ be a gamma function on the group $(G,\cdot)$ that takes
  values in $\Inn(G)$.

  Can  the  skew  brace   $(G,\cdot,  \circ)$,  where  $g\circ  h=g\cdot
  \lexp{h}{\gamma(g)}$  for   all  $g,h\in   G$,  be  obtained   from  a
  Rota--Baxter operator $B\colon G\to G$?

  Equivalently,  
  does there exist a Rota--Baxter operator $B$ on $G$ such that
  \begin{equation}
    \label{RB-to-gamma}
    \gamma(g) = \iota(B(g))
  \end{equation}
  holds for all $g \in G$?
  
  In this case, we say that the gamma function \emph{comes} from a
  Rota--Baxter operator. 
\end{question}

\begin{remark}
  It might be noted that it is not difficult to find skew braces
  for which the image of the corresponding gamma function is
  \emph{not} contained in the group of inner automorphisms.

  For  instance,  this  is  the   case  when  $(G,\cdot,\circ)$  is  a
  non-trivial  skew   brace  with  $(G,  \cdot)$   abelian  (that  is,
  $(G,\cdot,\circ)$ is a non-trivial \emph{brace}; see~\cite{braces}).
  
  In Section~\ref{sec:nonabelian},
  we exhibit a class of examples where $(G, \cdot)$ is non-abelian.
\end{remark}

Question~\ref{q:titq} has  a positive answer  when $G$ is  a centreless
group. In fact, in this case the  map $\iota$ is an isomorphism, so if
$\gamma : G \to \Inn(G)$ is  a gamma function, then the composition of
$\gamma$  with $\iota^{-1}$  yields  a Rota--Baxter  operator $B$  for
which~\eqref{RB-to-gamma}   holds.   (Compare   with~\cite[Proposition
  3.13]{BarGub2}.)

In this paper,  we show  that  the question
has a  negative answer  in 
general, by exhibiting two  counterexamples, which are finite
$p$-groups, for $p$ an odd prime. 

Our examples  rely on  a cohomological  characterisation of  the gamma
functions that  come from  a Rota--Baxter operators,  in terms  of the
vanishing of a  certain cocycle in a suitable  second cohomology group
over   a   trivial   module.   Our  characterisation   is   given   in
Theorem~\ref{thm:thm},  and   depends  in   turn  on   a  well-known
cohomological  characterisation  (Proposition~\ref{prop:MOF})  of  the
morphisms from a group $U$ to a quotient of another group $V$ by an
abelian normal subgroup, that
can be lifted to a morphism $U \to V$.

In~Section~\ref{sec:extensions},  we  recall  the  standard  connection
between  group  extensions  and   the  second  cohomology  group.   In
Section~\ref{sec:reconstructing}, we  show that,  if a  gamma function
comes   from  a   Rota--Baxter   operator,   then  our   cohomological
characterisation allows us to reconstruct  the latter from the former,
on the basis  of the knowledge that a certain  central group extension
splits.

Our first example, given in Section~\ref{sec:example}, 
has order $p^{5}$,  and it is a simplified
version  of   an  example  of~\cite{CS2}. With  the  cohomological
setting  in  place, the  proof that  this group provides a  negative
answer to  Question~\ref{q:titq} reduces  to the elementary  fact that
the Heisenberg group of order $p^{3}$ does not split over its centre.

In Section~\ref{sec:more}, we study a family of gamma functions on the
Heisenberg group  of order  $p^{3}$, and  show,  once more  via the
connection with  the splitting of certain extensions,  that precisely one
of them provides another, smaller counterexample
to  Question~\ref{q:titq}.   For the  other  gamma  functions in  this
family,  we apply  the methods  of Section~\ref{sec:reconstructing}  to
find explicitly the corresponding Rota--Baxter operators.

\section{Group extensions and the second cohomology group}
\label{sec:extensions}

We recall here the connection  between group extensions and a suitable
second  cohomology  group;  see~\cite[Chapter IV,  Section  3]{Bro82},
\cite[Section 5.1]{Wei69}.

Let  $G$ be  a group,  let $Q$  be a  $G$-module, where  the
action is written as left exponent,  and let $E$ be an \emph{extension
  of $G$  by $Q$},  that is,  a group  $E$ together  with an
exact sequence of groups
\begin{equation}
  \label{eq:exact}
  1 \to Q \to E \xrightarrow{\pi} G\to 1
\end{equation}
such that the  action by conjugation of $G\cong E/Q$  on $Q$ coincides
with the original  one. (Here we identify $Q$ with  its image in $E$.)
Let  $s  \colon   G  \to  E$  be  a  section,   that  is,  a
set-theoretic map $G\to E$ such that $\pi(s(x)) = x$ for all
$x  \in G$.  Then  there exists  a  function $\theta  \colon
G \times  G \to  Q$ such  that for  all $a,  b \in
G$,
\begin{equation}
  \label{eq:theta}
  s(a) s(b) = \theta(a, b) s(a b) .
\end{equation}
Associativity now yields that $\theta\colon G \times
G \to Q$ is  a $2$-cocycle: for all $a,b,c\in G$,
\begin{equation*}
  \lexp{\theta(b,c)}{a}\theta(a,bc)=\theta(ab,c)\theta(a,b). 
\end{equation*}

Conversely, if $\theta\colon G \times
G \to Q$ is
a $2$-cocycle,
then one may endow
the set $E = Q \times G$ with the group operation
\begin{equation}
  \label{eq:op-E}
  (q_{1}, a_{1}) (q_{2}, a_{2})
  =
  (q_{1} (\lexp{q_{2}}{a_{1}}) \theta(a_{1}, a_{2}), a_{1} a_{2})
\end{equation}
and      obtain    an   exact    sequence~\eqref{eq:exact},
with the natural injection on the first component and
  projection on the second one, for which~\eqref{eq:theta} holds with
$s(a) = (1, a)$ for all $a\in G$.

Putting all these facts together, we have obtained most of the following
standard result. 
\begin{proposition}
  \label{prop:rob}
  Let $G$ be a group, and let $Q$ be a $G$-module. 
  \begin{enumerate}
  \item
    Given an extension $E$ of $G$  by $Q$ as in~\eqref{eq:exact}, and a
    section  $s$,  equation~\eqref{eq:theta}   defines  a  $2$-cocycle
    $\theta$, whose class  in $\coh^{2}(G, Q)$ does not  depend on the
    particular section we have chosen.
  \item
    If    $\theta\colon G \times  G \to  Q$ is  a
    $2$-cocycle, then there  exist an
    exact sequence~\eqref{eq:exact} and a section $s$ such
    that equation~\eqref{eq:theta} holds. 
  \item
    The following are equivalent:
    \begin{enumerate}
    \item 
      The extension $E$ of $G$ by $Q$ splits, that is, there
      exists a section which is a morphism.
    \item
      The cocycle $\theta$ yields the trivial class in
      $\coh^{2}(G, Q)$. 
    \end{enumerate}
  \end{enumerate}
\end{proposition}

\begin{remark}
  We will apply this result only  for trivial modules. In particular, if
  $Q$ is  a trivial $G$  module, then the corresponding  group extension
  $E$ is  \emph{central}, meaning that (the  image of) $Q$ is  contained in the
  centre of $E$.

  Moreover, note that if a central extension $E$ of $G$ by $Q$ splits,
  then $E\cong G\times Q$.
\end{remark}

\section{Rota--Baxter operators and cohomology}
\label{sec:cohomology}

Proposition~\ref{prop:MOF}  below   is  well known;  it   follows  for
instance from  the diagrammatic approach of~\cite[Chapter  IV, Section
  3,  Exercise 1(a)]{Bro82},  and  Proposition~\ref{prop:rob}.  In  the
following we sketch a direct, elementary argument.

Let $U$ and $V$  be groups, let $A$ be an  abelian, normal subgroup of
$V$, and let $\psi \colon U \to V/A$ be a morphism. Note that $A$ is a
$V/A$-module  under  conjugation, and  then  a  $U$-module via
$\psi$; we denote  this action by a left exponent.

Let $C\colon U \to V$ be a lifting of $\psi$,
that is, a set-theoretic map such that $\psi(u) = C(u) A$ for all $u \in
U$. Since $\psi$ is a morphism, we have for all $u_{1}, u_{2} \in U$,
\begin{equation}
  \label{eq:enter-coc}
  C(u_{1}) C(u_{2}) = \kappa(u_{1}, u_{2}) C(u_{1} u_{2}),
\end{equation}
for a suitable function $\kappa \colon U \times U \to A$. Expanding
\begin{equation*}
  C((u_{1} u_{2}) u_{3})
  =
  C(u_{1} (u_{2} u_{3})),
\end{equation*}
one   sees    that   $\kappa\colon   U\times    U   \to   A$    is   a
$2$-cocycle. Moreover, it is immediate to  see that if we set, for all
$u\in U$,
\begin{equation*}
  C'(u) = \sigma(u) C(u),
\end{equation*} for a
function $\sigma \colon U \to A$, then the cocycle $\kappa'$ associated to
$C'$ is given by
\begin{equation*}
  \kappa'(u_{1}, u_{2})
  =
  \kappa(u_{1}, u_{2}) \sigma(u_{1})
  (\lexp{\sigma(u_{2})}{u_1}) \sigma(u_{1} u_{2})^{-1},
\end{equation*}
so it differs from $\kappa$ by a $2$-coboundary. 
From this discussion, we obtain the following well-known result.
\begin{prop}
  \label{prop:MOF}
  Let $U$  and $V$ be groups, let $A$ be an abelian, normal subgroup
  of $V$, and let 
  $\psi \colon U \to V/A$ be a morphism. Let $C\colon U \to V$ be a
  lifting of $\psi$, 
  and let $\kappa$ be the map $U \times U \to A$ defined
  by~\eqref{eq:enter-coc}. Then the following hold: 
  \begin{enumerate}
  \item
    $\kappa$ is a $2$-cocycle, whose class in $\coh^{2}(U, A)$ does
    not depend on the choice of the lifting $C$.
  \item
    The following are equivalent:
    \begin{enumerate}
    \item
      There exists a morphism $\phi \colon U \to V$ which is a lifting
      of $\psi$. 
    \item
      The class of $\kappa$ in $\coh^{2}(U, A)$ is trivial.
    \end{enumerate}
  \item
    Two morphisms $\phi_{1}, \phi_{2}$ are liftings of the same $\psi$ if
    and only if there exists a $1$-cocycle $\zeta \colon U \to A$ such
    that, for all $u \in U$,
    \begin{equation}
      \label{eq:same-gamma}
      \phi_{2}(u) =  \zeta(u) \phi_{1}(u).
    \end{equation}
    
  \end{enumerate}
\end{prop}

Let now  $\gamma$ be a  gamma function on  a group $(G,  \cdot)$, whose
values are  inner automorphisms,  and write $Z(G)$  for the  centre of
$(G,\cdot)$.
For  all  $g,h\in  G$,  let  \begin{equation*}  g  \circ  h  =  g\cdot
  \lexp{h}{\gamma(g)}.
\end{equation*}
Then  $(G, \circ)$  is  a  group and  $\gamma  \colon  (G, \circ)  \to
\Inn(G)$  is   a  morphism.   Composing  $\gamma$  with   the  natural
isomorphism  $\Inn(G)  \to  (G,  \cdot)/Z(G)$, we  obtain  a  morphism
$\psi\colon (G, \circ) \to (G, \cdot)/Z(G)$.

Note that a lifting  of $\psi$ to a morphism $B  \colon (G, \circ) \to
(G, \cdot)$ is precisely a Rota--Baxter operator on $G$ from which
$\gamma$ comes, as for such a
$B$,
\begin{equation*}
  B(g \cdot \lexp{h}{\iota(B(g))}) = B(g \circ h) = B(g) \cdot B(h)
\end{equation*}
for all $g,h\in G$.  Thus we can specialise Proposition~\ref{prop:MOF}
to the following result, where $Z(G)$ is  a trivial $(G,\circ)$-module. 
\begin{theorem}
  \label{thm:thm}
  Let  $\gamma$ be  a gamma  function on  the group  $G$, which  takes
  values  in $\Inn(G)$.  Let $C  :  G \to  G$  be a  function such  that
  $\gamma(g) = \iota(C(g))$ for all $g  \in G$, and let $\kappa \colon
  G  \times G  \to Z(G)$  be  defined by
  \begin{equation}
    \label{eq:coc}
    C(g)\cdot C(h)
    =
    \kappa(g, h) \cdot  C(g \circ h)
  \end{equation}
  for all $g, h \in G$.
Then the
  following hold:
  \begin{enumerate}
  \item
    $\kappa$ is a 
    $2$-cocycle, whose cohomology class 
    in 
    \begin{equation*}
      \coh^{2}((G, \circ), Z(G))
    \end{equation*}
    does not  depend on the choice  of $C$.
  \item
    The following are equivalent:
    \begin{enumerate}
    \item
      The gamma function $\gamma$ comes from a Rota--Baxter operator.
    \item
      The cohomology class of $\kappa$ in $\coh^{2}((G, \circ),
      Z(G))$ is trivial. 
    \end{enumerate}
  \item
    Two Rota--Baxter operators $B_{1}, B_{2}$ yield the same gamma
    function if and only if there exists a morphism $\zeta : (G, \circ)
    \to Z(G)$ such that
    \begin{equation*}
      B_{2}(g) = \zeta(g) \cdot B_{1}(g) 
    \end{equation*}
    for all $g \in G$.
  \end{enumerate}
\end{theorem}
Note that since
$Z(G)$ is a trivial $(G, \circ)$-module, the morphisms $\zeta$ are 
precisely the $1$-cocycles.

\section{Reconstructing the Rota--Baxter operator}
\label{sec:reconstructing}

With the same notation of the previous section, we now show how one
can explicitly find the Rota--Baxter operator associated with the
gamma function, when the cohomology class of the corresponding
cocycle is trivial. 

Let   $(G,\cdot)$    be   a   group   with    centre   $Z(G)$,   let
$\gamma(g)=\iota(C(g))$ be  a gamma function  on $G$ with  values in
$\Inn(G)$,  corresponding  to  an  operation  ``$\circ$'',  and  let
$\kappa$ be the cocycle obtained as in~\eqref{eq:coc}.

Assume that the central extension 
\begin{equation*}
  1 \to Z(G)\to E \to (G,\circ)\to 1
\end{equation*}
associated with  $\kappa$ splits, that is, that there is  a section
$s'\colon (G,\circ)\to E$  which is a  morphism. The section
$s'$ yields  the trivial  cocycle, and  to write  $\kappa$, associated
with the section $s(g)=(1,g)$, as  a coboundary, consider the map
$\sigma\colon  G\to   Z(G)$  such  that  $s'(g)= \sigma(g) s(g)  $.  We
immediately deduce that for all $g,h\in G$,
\begin{equation*}
  \kappa(g,h)=\sigma(g)^{-1}\cdot \sigma(h)^{-1}\cdot\sigma(g\circ h).
\end{equation*}
If   we now define
$B(g)= \sigma(g) \cdot C(g)  $  for   all   $g\in  G$,   then  by an
immediate application of Theorem~\ref{thm:thm},  $B$ is  the 
Rota--Baxter operator on $G$ from which $\gamma$ comes.

\section{An example of order $p^{5}$}
\label{sec:example}

This is a simplified version  of~\cite[Example 5.4]{CS2}.

Let $p$ be an odd prime, and let $H$ be the Heisenberg group of order
$p^{3}$: 
\begin{equation*}
  H
  =
  \Span{ u, v, k : u^{p}, v^{p}, k^{p},
    [u, v] = k, [u, k], [v, k]};
\end{equation*}
every element of $H$ can be written uniquely as
\begin{equation*}
	u^i v^j k^{q},
\end{equation*}
with $0\le i,j,q<p$.

Let  $(G,\cdot)  = S \times  H$,  where $S = \Span{x,  y}$  is
elementary  abelian   of  order  $p^{2}$,   so  that  $G$   has  order
$p^{5}$.  Write $Z(G)=Z(G,\cdot)$ and   $K =  \Span{k}=Z(H) \le
Z(G)$. Since $H/K = \Span{ u K,  v K }$ is elementary abelian of order
$p^{2}$, the assignment
\begin{equation*}
 \psi (x K) = u K,
  \quad
  \psi(y K) = v K,
  \quad
  \psi(h K) = K
\end{equation*}
for all $h\in H$ uniquely defines an  endomorphism $\psi$ of $G/K$,
whose kernel and 
image are $H/K$.

Let $C : G \to G$ be any function such that 
\begin{equation}
  \label{eq:itsalift}
  C(g) K = \psi(g K)
\end{equation}
for all $g  \in  G$.   (Such  a  function  is  denoted  by  $\lift{\psi}$
in~\cite{CS2}.) Note that we do not assume $C$ to be constant on the
cosets of $K$. Also, note that $C(G) \subseteq H$ and $C(H)\subseteq K$. 

It was  proved in~\cite{CS2}, and it is immediate to see, that  the function
$\gamma : G \to \Inn(G)$  defined by $\gamma(g) = \iota(C(g))$ depends
only on $\psi$, and  not on the particular choice of  $C$.  It was also
proved in~\cite[Theorem 3.2]{CS2}  that $\gamma$  is a  gamma
function,  which thus 
defines a  skew brace  $(G, \cdot, \circ)$,  where ``$\cdot$''  is the
original  group  operation  on  $G$  and  ``$\circ$''  is  the  group
operation of~\eqref{eq:circ}, namely
\begin{equation*}
  a \circ b
  =
  a\cdot C(a)\cdot b\cdot C(a)^{-1}
\end{equation*}
for all $a,b\in G$.
(Actually, gamma functions are not mentioned explicitly in~\cite{CS2},
but they are used implicitly via~\eqref{eq:circ}.)

Note that  $[a, b] = 1$  for all $a \in S$
  and $b \in  C(G) \subseteq H$; therefore  the operations ``$\cdot$''
  and ``$\circ$'' coincide on $S$, which is thus a subgroup
  of $(G, \circ)$.

We now compute the cohomology class of the  $2$-cocycle $\kappa$
associated to $\gamma$, as per Section~\ref{sec:cohomology}. To
compute this class,  we are free to choose any of the functions $C$
satisfying~\eqref{eq:itsalift}. Our choice is the map
$C : G \to G$ defined by
\begin{equation*}
  C(x^{i} y^{j} c)
  =
  u^{i} v^{j}
  \in H
\end{equation*}
for all $0 \le i, j < p$ and $c \in H$. Let us compute the relevant
$2$-cocycle. We have, for all 
$0 \le i, j, m, n < p$ and $c, d \in H$,
\begin{align*}
C((x^{i} y^{j} c) \cdot  C(x^{i} y^{j} c) \cdot (x^{m} y^{n} d)
  \cdot C(x^{i} y^{j} c)^{-1})&=	 C(x^{i+m} y^{j+n} e)\\
  &= u^{i+m} v^{j+n},
\end{align*}
for some $e \in H$.
On the other hand,
\begin{align*}
  C(x^{i} y^{j} c) \cdot C(x^{m} y^{n} d)
  &=
  u^{i} v^{j} u^{m} v^{n}
  \\&=
  u^{i+m} v^{j+n} [v^{-j}, u^{-m}]
  \\&=
  u^{i+m} v^{j+n} k^{-j m}.
\end{align*}
So the relevant $2$-cocycle here is
\begin{equation}
  \label{eq:kappa}
  \kappa( x^{i} y^{j} c, x^{m} y^{n} d)
  =
  k^{-j m},
\end{equation}
with image in $K = Z(H) \le Z(G) = S \times K$.

The skew  brace $(G,  \cdot, \circ)$
provides a negative answer to  Question~\ref{q:titq}. This will follow
from Theorem~\ref{thm:thm} and the following lemma.
\begin{lemma}
  \label{lemma:non-trivial}
  The cocycle $\kappa$ of~\eqref{eq:kappa} yields a non-trivial class in
  the $2$-cohomology group $\coh^{2}( (G, \circ), Z(G))$.
\end{lemma}

We give two proofs of this lemma. Both rely on  the 
connection between central extensions  and a suitable second
cohomology group, as recalled in Section~\ref{sec:extensions}.

\begin{proof}[First proof of Lemma~\ref{lemma:non-trivial}]
  Consider $Z(G)=S\times K$ as a trivial module for $(G,
  \circ)$ and $K$ as a trivial module for $S$.
  
  As the actions are trivial, the pair of morphisms
  \begin{equation*}
    \iota\colon H\hookrightarrow G,\qquad \phi\colon S\times K\to K,
  \end{equation*}
  the inclusion and the projection respectively, yields a morphism
  \begin{equation*}
    \coh^{2}( (G, \circ), Z(G))
    \to
    \coh^{2}( (S, \circ), K);
  \end{equation*}
  see, for example,~\cite[Chapter III, section 8]{Bro82}.
  Explicitly, the cohomology class of the $2$-cocycle $\kappa$ is mapped to the
  cohomology class of the $2$-cocycle $\kappa'$ defined by
  \begin{equation}
    \label{eq:kappaprime}
    \kappa'
    (
    x^{i} y^{j},
    x^{m} y^{n}
    )
    =
    k^{-j m}
  \end{equation}
  for all $0\le i,j,m,n< p$. 

  The surjective morphism
  \begin{align*}
    \pi\colon H &\to S\\
    u^{i} v^{j} c & \mapsto x^{i} y^{j},
  \end{align*}
  with $0 \le i, j < p$ and $c\in K$, has kernel $K$,
  so it yields the central extension 
  \begin{equation*}
    1 \to K \hookrightarrow H \xrightarrow{\pi} S \to 1.
  \end{equation*} 
  Now  patently the  Heisenberg  group  $H$ does  not
    split over its centre $K$.
    
  We now show that the cohomology class associated  to this central
  extension  is  that   of  $\kappa'$. According  to
  Proposition~\ref{prop:rob},  this  will  yield that  $\kappa'$  is
  non-trivial in cohomology.  It will follow that the same holds for
  $\kappa$, thereby concluding the proof of the lemma.

  Consider the following section of $\pi$: 
  \begin{align*}
    s \colon  S&\to H\\
    x^{i} y^{j} &\mapsto u^{i} v^{j},
  \end{align*}
  where $0 \le i,j < p$. As
  \begin{align*}
    s(x^{i} y^{j}) s(x^{m} y^{n})= u^{i} v^{j} u^{m} v^{n}
    =
    u^{i+m} v^{j+n}  [u,v]^{-j m} = u^{i+m} v^{j+n}  k^{-j m}
  \end{align*}
  and 
  \begin{equation*}
    s((x^{i} y^{j}) (x^{m} y^{n} ))
    =
    s(x^{i+m} y^{j+n})=u^{i+m} v^{j+n}
  \end{equation*}
  for all $0\le i,j,m,n<p$, 
  the cocycle we are looking for is
  \begin{equation*}
    (x^{i} y^{j}, x^{m} y^{n}) \mapsto k^{-jm},
  \end{equation*}
  as claimed.
\end{proof}

We  now give  an alternative  proof of  Lemma~\ref{lemma:non-trivial},
which is  elementary, in  that it  avoids the  use of
the above map between cohomology groups, replacing it
with a group-theoretic argument.

\begin{proof}[Second proof of Lemma~\ref{lemma:non-trivial}]
  Let $E = (S \times K) \times (G, \circ)$, and consider the central
  extension  
  \begin{equation}
    \label{eq:alt}
    1 \to S \times K \to E \to (G, \circ) \to 1
  \end{equation}
  defined by $\kappa$  and the standard section $s(g) =  (1, g)$. Thus
  the operation  on $E$  is given  by~\eqref{eq:op-E}, with  $\theta =
  \kappa$; it  follows from this and  the formula~\eqref{eq:kappa} for
  $\kappa$ that
  \begin{align*}
    [(1, x), (1, y)]
    &=
    (1, x) (1, y) ((1, y) (1, x))^{-1}
    \\&=
    (1, x y) (\kappa(x, y), 1) ( (\kappa(y, x), 1) (1, y x) )^{-1}
    \\&=
     (1, x y) (1, y x)^{-1} (\kappa(y, x), 1)^{-1}
    \\&=
    (\kappa(y, x)^{-1}, 1)
    =
    (k, 1),
  \end{align*}
  so that $(\nothing{1} \times K, \nothing{1})$ is contained in the
  derived subgroup 
  $[E, E]$ of $E$.

  Assume  by way  of  contradiction  that the  sequence~\eqref{eq:alt}
  splits, and let $T$ be a complement  to $S \times K$ in $E$. Then $M
  = (S  \times \nothing{1}, \nothing{1})  T$ is a maximal  subgroup of
  the  finite  $p$-group $E$,  which  does  not contain  $(\nothing{1}
  \times K, \nothing{1})$, a non-trivial subgroup of $[E, E]$.

  This contradiction shows that~\eqref{eq:alt} does
  not split, so that, by Proposition~\ref{prop:rob}, $\kappa$ is
  non-trivial in $\coh^{2}((G, \circ), Z(G))$.
\end{proof}

\section{Examples of order $p^{3}$}
\label{sec:more}

Consider again the Heisenberg group of order $p^{3}$, where $p$ is an odd prime:
\begin{equation*}
  (H,\cdot)
  =
  \Span{
    u, v, k
    :
    u^{p}, v^{p}, k^{p}, [u, v] = k, [u, k], [v, k]
  }.
\end{equation*}
Write, as usual, $Z(H)$ for the centre of $(H,\cdot)$. 

For   all  $\alpha \in \ZZ / p \ZZ$,  consider   the  function   $\gamma(g)  =
\iota(g^{\alpha})$, which is showed in~\cite[Proposition 5.6]{CS2} to be a
gamma function on $H$. The associated circle operation is
\begin{equation*}
  g \circ h
  =
  g\cdot \lexp{h}{\iota(g^{\alpha})}
  =
  g\cdot  h\cdot [g, h]^{\alpha}
\end{equation*}
for all $g,h \in H$. 
Note that powers with respect to $\circ$ (in particular, inverses)
coincide with the corresponding powers in $H$. Moreover, the commutator
in $(H, \circ)$ of $g$ and $h$ 
is given by
\begin{align*}
  [g,h]_{\circ}
  &=
  g\circ h\circ g^{-1}\circ h^{-1}=(g\cdot h\cdot [g,
    h]^{\alpha})\circ(g^{-1}\cdot h^{-1}\cdot
  [g^{-1},h^{-1}]^{\alpha})\\ 
  &= g\cdot h\cdot [g,h]^{\alpha}\cdot g^{-1}\cdot h^{-1}\cdot
  [g^{-1}, h^{-1}]^{\alpha}\cdot [g\cdot h,g^{-1}\cdot
    h^{-1}]^{\alpha}\\ 
  &=[g, h]^{1+2\alpha}\cdot[h\cdot g, g\cdot h]^{\alpha}
  =
  [g, h]^{1+2\alpha}
  .
\end{align*}
In particular, $(H, \circ)$ is abelian precisely when $\alpha = -1/2$,
where the fraction is taken in $\ZZ  / p \ZZ$. (Note that this situation
is  a particular  case  of  the Baer  correspondence~\cite{Baer-corr},
which is in turn an approximation of the Lazard correspondence and the
Baker--Campbell--Hausdorff formulas~\cite[Chapters~9
  and~10]{khukhro}.)

We now show the following result.
\begin{prop}
  \ppar
  \begin{enumerate}
  \item
    When $\alpha = -1/2$, the  gamma function $\gamma$ does \emph{not}
    come from a Rota--Baxter operator.
  \item
    When $\alpha \ne  -1/2$, the gamma function $\gamma$  comes from the
    Rota--Baxter operator
    \begin{equation}
      \label{eq:R-B-operator}
      B(u^{i}\cdot  v^{j}\cdot   k^{r})
      =
      u^{i\alpha}\cdot  v^{j\alpha}
      \cdot 
      k^{\alpha^{2} (r-i j \alpha)(1 + 2 \alpha)^{-1}},
    \end{equation}
    for $0 \le i,j,r < p$.
  \end{enumerate}
\end{prop}

When $\alpha =  -1/2$, we have  thus obtained another example  giving a
negative answer to Question~\ref{q:titq}.

\begin{proof}
Take $C \colon H \to H$ to be $C(g) = g^{\alpha}$. We have
\begin{align*}
  C(g \circ h)
  &=
  C(g\cdot  h\cdot  [g, h]^{\alpha})
  =
  g^{\alpha}\cdot h^{\alpha} \cdot [h, g]^{\binom{\alpha}{2}}\cdot  [g, h]^{\alpha^2}
  \\&=
  C(g)\cdot  C(h)\cdot  [g, h]^{\frac{\alpha^2+\alpha}{2}}.
\end{align*}
The critical cocycle $\kappa \colon H \times H \to Z(H)$ is thus
\begin{align*}
  \kappa(g, h)
  =
  [g, h]^{-\frac{\alpha^2+\alpha}{2}}
\end{align*}
for all $g,h\in H$.
Consider, as in Section~\ref{sec:extensions}, the standard central
extension defined by $\kappa$: 
\begin{equation}
  \label{eq:ses}
  \begin{aligned}
    1 \to Z(H) \to Z(H) \times H \to (H, \circ) \to 1,
  \end{aligned}
\end{equation}
where the operation on $E = Z(H) \times H$ is given by
\begin{align*}
  (a, g) (b, h)
  =
  (a\cdot b\cdot \kappa(g, h), g \circ h)
\end{align*}
for all $a,b\in Z(H)$ and $g,h\in H$, with the standard section $s : H
\to E$ given by $s(g) = (1, g)$ for all $g\in H$. 

Let us compute the commutator of two elements $(1, g), (1, h)$ in
$E$. We have
\begin{align*}
  [(1, g), (1, h)]
  &=
  (1, g) (1, h)((1, h) (1, g))^{-1}
  \\&=
  (\kappa(g, h), g \circ h)
  (\kappa(h, g), h \circ g)^{-1}
  \\&=
  (\kappa(g, h), g \circ h)
  (\kappa(h, g)^{-1}, (h \circ g)^{-1})
  \\&=
  (\kappa(g, h)^{2}, [g, h]_{\circ})
  \\&=
  ([g, h]^{-\alpha (\alpha + 1)}, [g, h]^{1 + 2 \alpha}),
\end{align*}
where
we have
exploited
\begin{equation*}
  \kappa((g \circ h),(h \circ g)^{-1})
  =
  \kappa( g\cdot  h, (h \cdot g)^{-1})
  =
  1.
\end{equation*}
It is  immediate to see that  $(Z(H), Z(H)) \le Z(E)$.  (The following
arguments yield that they  actually coincide.) Therefore the commutator
subgroup $[E, E]$ is of order $p$, generated by
\begin{equation}
  \label{eq:comm-uv}
  [(1, u), (1, v)]
  =
  (k^{-\alpha (\alpha + 1)}, k^{1 + 2 \alpha}).
\end{equation}
When $\alpha =  -1/2$, we have seen above that  the group $(H, \circ)$
is abelian, so~\eqref{eq:ses} does not  split, as otherwise $E \cong Z(H)
\times (H, \circ)$ would be abelian, whereas~\eqref{eq:comm-uv} shows that
$[E, E]$ is non-trivial.

Consider now the case $\alpha \ne -1/2$.  The previous argument yields
that both  $(H,\circ)$ and  $S=\Span{(1, u),  (1, v)}$  are Heisenberg
groups of order $p^{3}$. The former  is generated by $u,v$, with $\widetilde
k=[u,v]_{\circ}=k^{1+2\alpha}$,  and   the  latter  is   generated  by
$(1,u),(1,v)$, with $[(1,u),(1,v)]=(k^{-\alpha (\alpha + 1)}, k^{1 + 2
  \alpha})$.
An isomorphism between $(H,\circ)$ and $S$ is clearly given by 
\begin{align*}
  s'\colon(u^i\circ v^j\circ \widetilde{k}^q)
  \mapsto
  &\ (1,u)^i(1,v)^j(k^{-\alpha (\alpha + 1)}, k^{1 + 2 \alpha})^q
  \\&=(k^{-\frac{\alpha^2 + \alpha}{2}(2q+ij)},u^{i}\circ v^{j}\circ \widetilde{k}^q),
\end{align*}
with $0\le i,j,q<p$.  This implies immediately that $s'$  is a section
which splits the extension~\eqref{eq:ses}.

We now  employ the methods of  Section~\ref{sec:reconstructing} to compute
the Rota--Baxter operator $B$ from which $\gamma$ comes.
As
\begin{align*}
  s'(u^i\circ v^j\circ \widetilde{k}^q)&=(k^{-\frac{\alpha^2 +
      \alpha}{2}(2q+ij)},u^{i}\circ v^{j}\circ
  \widetilde{k}^q)\\ 
  &=(k^{-\frac{\alpha^2 + \alpha}{2}(2q+ij)},1)(1,u^{i}\circ
  v^{j}\circ \widetilde{k}^q)\\ 
  &=(k^{-\frac{\alpha^2 + \alpha}{2}(2q+ij)},1)
  s(u^{i}\circ v^{j}\circ \widetilde{k}^q),
\end{align*}
we have
\begin{equation*}
  B(g)=\sigma(g)\cdot C(g), 
\end{equation*}
where 
\begin{equation*}
  \sigma(u^{i}\circ v^{j}\circ \widetilde{k}^q)=k^{-\frac{\alpha^2 + \alpha}{2}(2q+ij)}
\end{equation*}
and 
\begin{equation*}
  C(u^{i}\circ v^{j}\circ \widetilde{k}^q)=(u^i\cdot v^j\cdot k^{ij\alpha+q(1+2\alpha)})^{\alpha}=u^{i\alpha}\cdot v^{j\alpha}\cdot k^{q(\alpha+ 2\alpha^2)+ij(\frac{\alpha^2+\alpha}{2})}.
\end{equation*}
Explicitly,
\begin{align*}
  B(u^i \circ v^j \circ \widetilde{k}^q)&=
  k^{-\frac{\alpha^2 + \alpha}{2}(2q+ij)}\cdot u^{i\alpha}\cdot v^{j\alpha}\cdot k^{q(\alpha+ 2\alpha^2)+ij(\frac{\alpha^2+\alpha}{2})}\\
  &= u^{i\alpha} \cdot v^{j\alpha}\cdot k^{q\alpha^2}.
\end{align*}
The expression~\eqref{eq:R-B-operator} for $B$ in terms of the
operation ``$\cdot$'' is then
easily obtained (see Remark~\ref{rem:the_switch} below).

We can double-check that $B$ is a Rota--Baxter operator:
\begin{align*}
  B(u^i \circ v^j \circ \widetilde{k}^q\circ u^m \circ v^n \circ \widetilde{k}^r)&=B(u^{i+m} \circ v^{j+n} \circ \widetilde{k}^{q+r-jm})\\
  &= u^{(i+m)\alpha} \cdot v^{(j+n)\alpha} \cdot k^{(q+r-jm)\alpha^2}
\end{align*}
and
\begin{align*}
  B(u^i \circ v^j \circ \widetilde{k}^q)\cdot B(u^m \circ v^n \circ
  \widetilde{k}^r)&=u^{i\alpha} \cdot v^{j\alpha}\cdot
  k^{q\alpha^2}\cdot u^{m\alpha} \cdot v^{n\alpha}\cdot
  k^{r\alpha^2}\\ 
  &=u^{i\alpha+ m\alpha}\cdot v^{j\alpha+ n\alpha}\cdot
  k^{-jm\alpha^2}\cdot k^{q\alpha^2+ r\alpha^2}, 
\end{align*}
so the assertion follows. 
\end{proof}

\begin{remark}
  \label{rem:the_switch}
  For $\alpha \ne  -1/2$, it is  easy to switch  from a writing  of the
  kind $u^{m}\cdot v^{n}\cdot k^{r}$ to one of the kind $u^i \circ v^j
  \circ \widetilde{k}^q$. Indeed,
  \begin{equation*}
    u^i \circ v^j \circ \widetilde{k}^q
    =
    u^i\cdot v^j\cdot k^{ij\alpha}\cdot k^{q(1+2\alpha)}
    =
    u^i\cdot v^j\cdot k^{ij\alpha+q(1+2\alpha)},
  \end{equation*}
  so that
  \begin{equation*}
    u^{i}\cdot v^{j}\cdot k^{r}
    =
    u^i \circ v^j \circ \widetilde{k}^q,
  \end{equation*}
  where $q = (r-i j \alpha)(1 + 2 \alpha)^{-1}$.
\end{remark}

\section{A class of gamma functions whose values are not all inner}
\label{sec:nonabelian}

Let $A$ and $B$ be two groups,  and let $\psi \colon A \to \Aut(B)$ be
a morphism such that $\psi(A) \nsubseteq \Inn(B)$. Write the action of
$\psi(A)$ on $B$  as left exponent, and consider  two group structures
on the set $G = A \times B$:
\begin{enumerate}
\item $(G, \cdot)$ is the direct product.
\item $(G, \circ)$ is the semidirect product $A \ltimes_{\psi} B$.
\end{enumerate}
Then  $(G,   \cdot,  \circ)$   is  a  skew   brace;  see~\cite[Example
  1.4]{skew}. Here the gamma function is given by
\begin{align*}
  \lexpp{a', b'}{\gamma(a , b)}
  &=
  (a , b)^{-1}\cdot ((a , b)\circ(a' , b'))
  \\&=
  (a , b)^{-1} \cdot (a \cdot a' , b \cdot \lexp{b'}{\psi(a)})
  \\&=
  (a' , \lexp{b'}{\psi(a)}),
\end{align*}
for all $a, a'\in A$ and $b,  b'\in B$.  Clearly $\gamma(a , b) \in
\Inn(G, \cdot)$ if and only if $\psi(a) \in \Inn(B)$. Moreover, $(G,
\cdot)$ is non-abelian if one of $A$ and $B$ is non-abelian.

To specify a concrete example, let  $A = \Span{a}$ be the cyclic group
of  order $4$,  and let $B  = D_4  =\Span{r,s: r^4=s^2=1,srs=r^{-1}}$ be  the
dihedral group of order  $8$. If $\psi(a)$ maps $r$ in  $r$ and $s$ in
$rs$,   then  $\psi(a)\notin   \Inn(B)$,  so   that  $\gamma(a,b)\notin
\Inn(A\times B,\cdot)$ for all $b\in B$.

As another  example, take  $V$ to  be an  elementary abelian  group of
order $p^{2}$ and $H  = \Span{h_{0}}$ to be a cyclic  group of order a
prime $q$, with $q \mid p - 1$. Choose two elements $\mu$ and $\nu$ of
order $q$ in $\Aut(V)\cong \GL(2, p)$  such that $\mu$ is scalar while
$\nu$  is not,  so that  $\mu$  and $\nu$  commute, but  they are  not
conjugate in $\GL(2, p)$.  Define $B$  to be the semidirect product of
$V$ by $H$ via $\phi(h_{0}) = \mu$.  Now let $A = \Span{a}$ be another
group  of order  $q$, and  take  in the  construction above  $\lexpp{h
  ,v}{\psi(a)} = (h , \nu(v))$ for all $h \in H$ and $v \in V$.

\providecommand{\bysame}{\leavevmode\hbox to3em{\hrulefill}\thinspace}
\providecommand{\MR}{\relax\ifhmode\unskip\space\fi MR }
\providecommand{\MRhref}[2]{%
  \href{http://www.ams.org/mathscinet-getitem?mr=#1}{#2}
}
\providecommand{\href}[2]{#2}

\end{document}